\documentclass[11pt]{article}
\makeatletter

\@addtoreset{equation}{section}
\makeatother
\usepackage{amsfonts}
\usepackage{eucal}
\usepackage{amsmath}
\usepackage{amssymb}
\newtheorem{theorem}{\bf Theorem}[section]
\newtheorem{definition}[theorem]{\bf Definition}

\newtheorem{lemma}[theorem]{\bf Lemma}

\newtheorem{remark}[theorem]{\bf Remark}
\newtheorem{example}[theorem]{\bf Example}
\begin{document}
\title{\bf 
Brudno's theorem for $\mathbb{Z}^d$ (or $\mathbb{Z}^d_+$) subshifts
}
  \author{Toru Fuda\footnote{Corresponding author. E-mail:
t-fuda@math.sci.hokudai.ac.jp} 
and Miho Tonozaki\footnote{Present affiliation: NEC Corporation, 
5-7-2, Shiba, Minato, Tokyo, Japan.}\\
Department of Mathematics\\
Hokkaido University\\
Sapporo 060-0810\\
Japan
}
\date{}

\maketitle

\begin{abstract}
We generalize Brudno's theorem of $1$-dimensional shift dynamical system 
to $\mathbb{Z}^d$ (or $\mathbb{Z}_+^d$) subshifts. 
That is to say, in $\mathbb{Z}^d$ (or $\mathbb{Z}^d_+$) 
subshift, the Kolmogorov-Sinai entropy is 
equivalent to the Kolmogorov complexity density almost everywhere 
for an ergodic shift-invariant measure.
\end{abstract}

\medskip

\begin{flushleft}
{\bf Mathematics Subject Classification (2010).} 
68Q30, %Algorithmic information theory (Kolmogorov complexity, etc.) 
28Dxx, %Measure-theoretic ergodic theory
37B10, %Symbolic dynamics
82B20. %Lattice systems (Ising, dimer, Potts, etc.) and systems on graphs
\\
\medskip
{\bf Keywords.} Brudno's theorem, Kolmogorov-Sinai entropy, Kolmogorov complexity, subshifts, 
$\mathbb{Z}^d$-action, pressure.

\end{flushleft}

\section{Introduction}

In a topological dynamical system, 
A. A. Brudno defined 
a complexity of the trajectory of a point in the space 
by using the notion of Kolmogorov complexity, 
and showed the relationship between this quantity and 
the Kolmogorov-Sinai entropy \cite{Br}.
As a preliminary step, Brudno 
considered the $1$-dimensional shift dynamical system and showed that, 
for an ergodic shift-invariant measure, 
the Kolmogorov complexity density is equal to the 
Kolmogorov-Sinai entropy almost everywhere \cite[Theorem 1.1]{Br}. 
%In what follows, we concentrate exclusively on 
%shift dynamical Brudno's theorem.

A partial approach to generalize this theorem to a $d$-dimensional 
case is found in \cite{Sim}.
S. G. Simpson showed that, in 
$\mathbb{Z}^d$ (or $\mathbb{Z}_+^d$) subshifts, there exists a point such 
that its Kolmogorov complexity density is coincident with the 
topological entropy \cite{Sim}. 
Examining Simpson's proof, we see that what
he showed substantively is that 
the Kolmogorov complexity density is equal to the 
Kolmogorov-Sinai entropy almost everywhere 
only for a measure of maximal entropy.

The purpose of this paper is to generalize the Brudno's theorem of 
the $\mathbb{Z}_+^1$-action
 shift dynamical system to $\mathbb{Z}^d$ (or 
$\mathbb{Z}_+^d$) subshifts. The main theorem is the following:

\setcounter{section}{3}
\begin{theorem}[Brudno's theorem for $\mathbb{Z}^d$ (or 
	$\mathbb{Z}^d_+$) subshifts]%\label{main}
	If $\mu \in EM(S, \varsigma)$, then 
	\begin{equation}
		%\lim_{n\to \infty}
		%\frac{\mathsf{K}(\omega \upharpoonright \Lambda_n)}{|\Lambda_n|}
		\mathcal{K}(\omega)=h_{\varsigma}(\mu), \quad
		\mu\text{-a.e.} \omega \in S.
	\end{equation}\label{Brudno}
\end{theorem}
\setcounter{section}{1}

Here $S$ denotes $\mathbb{Z}^d$ (or $\mathbb{Z}^d_+$) subshift, $\varsigma$ 
denotes the shift action on $S$, $EM(S,\varsigma)$ denotes the set of all ergodic 
shift-invariant measures on the topological dynamical system $(S,\varsigma)$, 
$\mathcal{K}(\omega)$ denotes the Kolmogorov complexity density of $\omega$, 
and $h_{\varsigma}(\mu)$ denotes the Kolmogorov-Sinai entropy of 
the measure preserving dynamical system 
$(S, \mathfrak{B}(S), \mu, \varsigma)$. 
We give the rigorous definition of these terms in Section 2.

In Section 2, we introduce some basic mathematical notions in 
ergodic theory, Kolmogorov complexity and shift dynamical systems.
We used \cite{Ke} and \cite{LiVi} as main references for this section.
Using these basic notions, we 
define the Kolmogorov complexity density 
of each point of $\Sigma^{\mathbb{Z}^d}$ 
(or $\Sigma^{\mathbb{Z}^d_+}$) naturally.

In Section 3, we prove the main theorem. 
The proof directly uses an idea of Brudno's original paper, i.e., 
Shannon-McMillan-Breimann theorem and the notion of frequency set. 

In the last section, as an application of the main theorem, we show a 
variational principle using the Kolmogorov complexity density.

%%%%%%%%%%%%%%%%%%%%%%%%

\section{Some Mathematical Preliminaries}

We first give quick reviews for some mathematical results related to 
the main theorem.
Descriptions of this section are restricted to a minimum and 
all the contents in this section are well known.
We write 
$$
	\mathbb{N}=\{1,2,\cdots\},\ 
	\mathbb{Z}=\{\cdots, -2, -1, 0, 1, 2, \cdots\},\ 
	\mathbb{Z}_+=\{0,1,2,\cdots\}.
$$
For an arbitrary fixed $d\in \mathbb{N}$, we set
$G:=\mathbb{Z}^d$ or $G:=\mathbb{Z}_+^d$.
For all $n\in \mathbb{N}$, let
$$
	\Lambda_n:=\{g=(g_i)_{i=1}^d \in G : \forall i\in \{1,\cdots,d\}, 
	|g_i|<n\}.
$$
Then we have
$$
	|\Lambda_n|=\begin{cases}
					(2n-1)^d & (G=\mathbb{Z}^d),\\
					n^d & (G=\mathbb{Z}_+^d),
				\end{cases}
$$
where we denote by $|A|$ the cardinality of a set $A$.

\subsection{Ergodic theory}

\subsubsection{Measure preserving dynamical system}

First, let us define Kolmogorov-Sinai entropy, also known as 
measure-theoretic entropy.

\begin{definition}[Measure preserving dynamical system]
Let $(X,\mathfrak{B},\mu)$ be a probability space and 
$\mathcal{T}=(T^g)_{g\in G}$ be a family of maps on $X$ such that
	\begin{enumerate}
		\item $\mu$ is $\mathcal{T}$-invariant, i.e., 
		$\forall g\in G, \forall A\in \mathfrak{B}, 
		\mu(T^{-g}A)=\mu(A)$ (where $T^{-g}:=(T^g)^{-1}$); 
		\item $\mathcal{T}$ is a measurable 
		action of $G$ on $X$, i.e., $T^g:X\to X$ is measurable for all 
		$g\in G$, 
		$T^0=I_X$ (the identity map on $X$) and 
		$\forall g, g'\in G, T^{g+g'}=T^g\circ T^{g'}$
		(if $G=\mathbb{Z}^d$, then $(T^g)^{-1}=T^{-g}$ holds).
	\end{enumerate}
	We call such a quadruple $(X,\mathfrak{B},\mu,\mathcal{T})$ 
	a measure preserving dynamical system (m.p.d.s.).
	\end{definition}

%In what follow, unless otherwise noted let
%$(X,\mathfrak{B},\mu,\mathcal{T})$ be a m.p.d.s..

\begin{definition}
	Let $(X,\mathfrak{B},\mu,\mathcal{T})$ be a m.p.d.s..
	A $\mathfrak{B}$-measurable function $f$ on $X$ is said to be 
	$\mathcal{T}$-invariant$\mod \mu$ if and only 
	if $\forall g\in G, f\circ T^g=f$ ($\mu$-a.s.).
	A set $A\in \mathfrak{B}$ is said to be 
	$\mathcal{T}$-invariant$\mod \mu$ 
	if and only if $1_A$ is $\mathcal{T}$-invariant$\mod \mu$, where we 
	denote the characteristic function of $A$ by $1_A$.
	We write $\mathcal{I}_{\mu}(\mathcal{T}):=\{A\in \mathfrak{B} : 
	A \ \text{is} \ \mathcal{T}\text{-invariant}\mod \mu\}
	=\{A\in \mathfrak{B} : \forall g\in G, \mu(T^{-g}A\bigtriangleup A)=0\}$, 
	where $\bigtriangleup$ denotes the symmetric difference.
\end{definition}

\begin{theorem}[Birkhoff's ergodic theorem]
	Let $(X,\mathfrak{B},\mu,\mathcal{T})$ be a m.p.d.s..
	Then, for all $f\in L^1(X, \mu)$, there exists the limit
	$$
		\bar{f}(x):=\lim_{n\to \infty}\frac{1}{|\Lambda_n|}
		\sum_{g\in\Lambda_n}f(T^gx), \quad
		\mu\text{-a.s.} \ x,
	$$
	and $\bar{f}\in L^1(X,\mu)$. 
	Moreover, $\bar{f}$ is $\mathcal{T}$-invariant$\mod \mu$ and
	$$
		\forall A\in \mathcal{I}_{\mu}(\mathcal{T}), \
		\int_A\bar{f}d\mu = \int_Afd\mu.
	$$	
\end{theorem}

{\it Proof}. 
See \cite{Ke}.
\hfill $\square$

\begin{definition}[Ergodicity]
	Let $(X,\mathfrak{B},\mu,\mathcal{T})$ be a m.p.d.s.
	If for all $f\in L^1(X,\mu)$
	$$
		\bar{f}=\int_X fd\mu, \quad \mu\text{-a.s.} \ x
	$$
	holds, then the 
	m.p.d.s. $(X, \mathfrak{B}, \mu, \mathcal{T})$ 
	is said to be ergodic.
	In this case, $\mu$ is called an ergodic $\mathcal{T}$-invariant 
	probability measure on the measurable space $(X, \mathfrak{B})$.
\end{definition}

Although there are several equivalent conditions of ergodicity, 
only the above-mentioned condition is used in this paper.

\begin{definition}[$\mu$-partition]
	Let $(X, \mathfrak{B}, \mu)$ be a probability space.
	A family of measurable sets $\alpha =\{A_i : i\in I\}\subset 
	\mathfrak{B}$ is called 
	a $\mu$-partition of $X$ if the following conditions hold:
	$$
		\mu(A_i\cap A_j)=0 \ (i\neq j), \ 
		\mu\left(X\setminus \bigcup_{i\in I}A_i\right)=0 \ 
		\text{and} \ 
		\mu(A_i)>0 \ (\forall i \in I).
	$$
	Accordingly, $\alpha$ is at most countable. 
	If $|I|<\infty$ is holds, then $\alpha$ is called a 
	finite $\mu$-partition.
\end{definition}

	Let $\alpha$ and $\beta$ be $\mu$-partitions of $X$. 
	The common refinement of $\alpha$ and $\beta$ 
	$$
		\alpha \vee \beta :=
		\{A \cap B : A\in \alpha , B \in \beta, \mu(A\cap B)>0\}
	$$
	is a $\mu$-partition of $X$.

\begin{definition}[Information and entropy of a $\mu$-partion]
	Let $(X, \mathfrak{B}, \mu)$ be a probability space, and $\alpha$ be a 
	$\mu$-partition of $X$. 
	The information of $\alpha$ is the measurable function $I_{\alpha}$ 
	on $X$ defined by 
	$$
		I_{\alpha}(x):=-\sum_{A\in\alpha}\log_2\mu(A)\cdot 1_A(x), 
		\quad x\in X.
	$$
	The entropy of $\alpha$ is defined by the average information, i.e.,
	$$
		H(\alpha):=\int_XI(\alpha)d\mu =\sum_{A\in \alpha}\varphi(\mu(A)),
	$$
	where we define the function $\varphi:[0,\infty)\to \mathbb{R}$ by
	$$
		\varphi(t):=	\begin{cases}
							-t\log_2t & (t>0),\\
							0 & (t=0).
						\end{cases}
	$$
\end{definition}

From Kolmogorov complexity's point of view, we choose the binary logarithm 
$\log_2$ instead of $\log_e$.

\begin{definition}[Dynamical entropy relative to a partition]
	Let $(X, \mathfrak{B}, \mu, \mathcal{T})$ be a m.p.d.s. and 
	$\alpha$ be a $\mu$-partition of $X$. We set 
	$T^{-g}\alpha :=\{T^{-g}A : A\in \alpha\}$ for each $g\in G$ and 
	$\alpha^{\Lambda}:=\bigvee_{g\in \Lambda}T^{-g}\alpha$ for a finite subset 
	$\Lambda \subset G$. 
	The dynamical entropy of the m.p.d.s. 
	$(X, \mathfrak{B}, \mu, \mathcal{T})$ relative to the partition $\alpha$ 
	is defined by
	$$
		h(\mu, \alpha, \mathcal{T}):=\inf_{n>0}\frac{1}{|\Lambda_n|}
		H(\alpha^{\Lambda_n}).
	$$
\end{definition}

\begin{theorem}
	Let $(X, \mathfrak{B}, \mu, \mathcal{T})$ be a m.p.d.s. and 
	$\alpha$ be a $\mu$-partition of $X$. Then
	$$
		h(\mu, \alpha, \mathcal{T})=\lim_{n\to \infty}\frac{1}{|\Lambda_n|}
		H(\alpha^{\Lambda_n}).
	$$
\end{theorem}
{\it Proof}. 
See \cite{Ke}.
\hfill $\square$

\begin{theorem}[Shannon-McMillan-Breiman]\label{SMB}
	Let $(X, \mathfrak{B}, \mu, \mathcal{T})$ be an ergodic m.p.d.s. and 
	$\alpha$ be a $\mu$-partition of $X$ with $H(\alpha)<\infty$. Then
	$$
		h(\mu, \alpha, \mathcal{T})=\lim_{n\to \infty}\frac{1}{|\Lambda_n|}
		I_{\alpha^{\Lambda_n}}
		\quad
		\text{in} \ L^1(X, \mu).
	$$
	Moreover, if $\alpha$ is finite, then this convergence 
	holds also for $\mu$-a.s.
\end{theorem}
{\it Proof}. 
See \cite{Ke, OW}.
\hfill $\square$

\begin{definition}[Kolmogorov-Sinai entropy]
	The Kolmogorov-Sinai entropy (KS entropy) of the m.p.d.s 
	$(X, \mathfrak{B}, \mu, \mathcal{T})$ is defined by
	$$
		h_{\mathcal{T}}(\mu):=
		\sup\{ h(\mu, \alpha, \mathcal{T}) : 
		\alpha \ \text{is a} \ \mu\text{-partition with} \ H(\alpha)<\infty\}.
	$$
\end{definition}

\begin{definition}[$\mu$-generator]
	Let $(X, \mathfrak{B}, \mu, \mathcal{T})$ be a m.p.d.s.. 
	A $\mu$-partition $\alpha$ is called a $\mu$-generator if 
	$\alpha^G=\mathfrak{B} \mod \mu$, where this equation means that 
	$\forall A\in \mathfrak{B}, \exists B \in \alpha^G, 
	\mu(A\bigtriangleup B)=0$.
\end{definition}

\begin{theorem}[Kolmogorov-Sinai]\label{KSth}
	Let $(X, \mathfrak{B}, \mu, \mathcal{T})$ be a m.p.d.s. 
	and $\alpha$ be a $\mu$-generator such that 
	$H(\alpha)<\infty$. Then $h_{\mathcal{T}}(\mu)
	=h(\mu, \alpha, \mathcal{T})$.
\end{theorem}
{\it Proof}. 
See \cite{Ke}.
\hfill $\square$

\subsubsection{Topological dynamical system}

We give the definition of topological dynamical system and its entropy via 
a variational principle of KS entropy.

\begin{definition}[Topological dynamical system]
	The pair $(X, \mathcal{T})$ 
	is called a topological dynamical system (t.d.s.) if the 
	following conditions hold:
	\begin{enumerate}
		\item $X$ is a compact metrizable space;
		\item $\mathcal{T}=(T^g)_{g\in G}$
		 is a continuous action of $G$ on $X$, i.e., 
		for all $g\in G$, $T^g:X\to X$ is continuous.
	\end{enumerate}
	In this setting we denote by $\mathfrak{B}(X)$ the Borel $\sigma$-algebra 
	of $X$. By (2), $\mathcal{T}$ is a measurable action of $G$ on $X$.
\end{definition}

Let $(X, \mathcal{T})$ be a t.d.s..
We denote by 
$M(X)$ the set of all probability measures 
on the Borel measurable space $(X, \mathfrak{B}(X))$, 
by $M(X, \mathcal{T})$ the set of all $\mathcal{T}$-invariant 
probability measures 
on $(X, \mathfrak{B}(X))$
and by $EM(X, \mathcal{T})$ the set 
of all ergodic members in $M(X, \mathcal{T})$, respectively.

\begin{theorem}[Krylov-Bogolubov]
	If $X\neq \emptyset$ and $(X, \mathcal{T})$ is a t.d.s., 
	then $M(X, \mathcal{T})\neq \emptyset$.
\end{theorem}
{\it Proof}. 
See \cite{Ke}.
\hfill $\square$

Obviously, if $\mu \in M(X, \mathcal{T})$ then $(X, \mathfrak{B}(X), \mu, 
\mathcal{T})$ is a m.p.d.s.. 
We define a concept of topological entropy of a t.d.s. $(X, \mathcal{T})$ 
as follows.

\begin{definition}[upper semicontinuous function]
	Let $Y$ be a topological space. 
		We set 
	$$
		USC(Y):=
		\left\{f:Y\to [-\infty, \infty): \forall c\in \mathbb{R}, \{y\in Y: f(y)<c\}\ \text{is open}
		\right\},
	$$
	and an element of $USC(Y)$ is called an upper semicontinuous function 
	on $Y$.
\end{definition}

\begin{definition}[Pressure, topological entropy, equilibrium state]
	Let $(X, \mathcal{T})$ be a t.d.s. and let $\psi \in USC(X), \ \inf \psi >-\infty$.
	The pressure of $\psi$ is defined by 
		$$
			p(\psi):=\sup_{\mu\in M(X, \mathcal{T})}(h_{\mathcal{T}}(\mu)+\mu(\psi))
		$$
	where $\mu(\psi):=\int_{X}\psi (x)d\mu(x)$. 
	A measure $\nu\in M(X,\mathcal{T})$ is called 
	an equilibrium state for $\psi \in USC(X)$ if 
		$$
			p(\psi)=h_{\mathcal{T}}(\nu)+\nu(\psi).
		$$
		In particular, $p(0)=\sup_{\mu\in M(X,\mathcal{T})}h_{\mathcal{T}}(\mu)$ is called the 
	topological entropy of $(X,\mathcal{T})$, and the equilibrium state for $\psi =0$ is called a measure of maximal 
	entropy for $\mathcal{T}$.
\end{definition}

\begin{theorem}[Ergodic decomposition]\label{ergdecomp}
	Let $(X, \mathcal{T})$ be a t.d.s.. Then, 
	for each $\mu \in M(X,\mathcal{T})$, there uniquely exists 
	a measure $\rho$ on the space 
	$M(X,\mathcal{T})$ (with respect to the Borel $\sigma$-algebra associated
	 to the weak-$\ast$ topology) such that
	\begin{enumerate}
		\item for any bounded measurable function $f:X\to \mathbb{R}$ we have
		$$
			\int_X f(x)d\mu(x)=\int_{EM(X,\mathcal{T})}
			\left\{
				\int_X f(x)d\nu(x)
			\right\}
			d\rho(\nu).
		$$
		\item $\rho(EM(X,\mathcal{T}))=1$.
	\end{enumerate}
\end{theorem}
{\it Proof}. 
See \cite{Ke, PY, Wal}.
\hfill $\square$

	Since 
	$\mu(A)=\int_{EM(X,\mathcal{T})}\nu(A)d\rho(\nu)$ for a measurable set $A\in \mathfrak{B}(X)$, 
	we write $\mu=\int_{EM(X,\mathcal{T})}\nu d\rho(\nu)$ which is called the ergodic decomposition 
	of $\mu$.

\begin{theorem}[Jacobs's theorem]\label{Jacob}
	Let $(X, \mathcal{T})$ be a t.d.s.. If $\mu \in M(X,\mathcal{T})$ and 
	$\mu=\int_{EM(X,\mathcal{T})}\nu d\rho(\nu)$ is the ergodic decomposition of $\mu$, then we have
	$$
		h_{\mathcal{T}}(\mu)=\int_{EM(X,\mathcal{T})}h_{\mathcal{T}}(\nu)d\rho (\nu).
	$$
\end{theorem}
{\it Proof}. 
See \cite{Ke, Wal}.
\hfill $\square$

%%%%%%%%%%%%%%%%%%%%%%%%%%%%%%%%%%

\subsection{Kolmogorov complexity}

Let $\mathcal{A}$ be a nonempty finite set. 
Without loss of generality, we set 
$\mathcal{A}:=\{0,1,\cdots, N\}$ where $N\in \mathbb{Z}_+$.

We define the set of all finite \textit{strings} 
over $\mathcal{A}$ as follows.
$$
	\mathcal{A}^{\ast}:=\bigcup_{n=0}^{\infty}\mathcal{A}^n
	=\{\lambda, 0,1,\cdots, N, 00,01,\cdots,0N, 10,\cdots ,1N, \cdots,
	NN, 000,\cdots \}
$$
where $\mathcal{A}^0=\{\lambda\}$ and $\lambda$ denote the empty string.

We often identify $\mathcal{A}^{\ast}$ with $\mathbb{Z}_+$ or $\mathbb{Z}$ 
by using the bijective map 
$I_{\mathcal{A}^{\ast}\to\sharp}:\mathcal{A}^{\ast}\to\sharp \ (\sharp \in \{\mathbb{Z}_+, \mathbb{Z}\})$ defined by the following.
\begin{eqnarray*}
	I_{\mathcal{A}^{\ast}\to\mathbb{Z}_+}(x)&:=&
	\begin{cases}
		\displaystyle\sum_{k=0}^{n-1}(N+1)^k+\sum_{k=1}^{n}a_k(N+1)^{n-k}, & x=a_1a_2\cdots
		a_n \in \mathcal{A}^n\ (n\in \mathbb{N}), \\
		0, & x=\lambda,
	\end{cases}\\
	I_{\mathcal{A}^{\ast}\to\mathbb{Z}}(x)&:=&
	\alpha\left(I_{\mathcal{A}^{\ast}\to\mathbb{Z}_+}(x)
	\right)
	\end{eqnarray*}
	where $\alpha(n):=(-1)^{n+1}\lfloor \frac{n+1}{2} \rfloor$ 
	for all $n\in \mathbb{Z}_+$.
	
	For example, the case of $\mathcal{A}=\{0,1\}$ is as follows:
	\begin{center}
	\begin{tabular}{c|cccccccccc}
		\hline
		$x$ & $\lambda$ & $0$ & $1$ & $00$ & $01$ 
		& $10$ & $11$ & $000$ & $001$ & $\cdots$
		\\ %\hline
		$I_{\{0,1\}^{\ast}\to\mathbb{Z}_+}(x)$ 
		& $0$ &$1$ & $2$ &$3$ & $4$
		& $5$ &$6$& $7$ & $8$ &  $\cdots$ \\
		\hline\hline
			$x$ & $\lambda$ & $0$ & $1$ & $00$ & $01$ 
		& $10$ & $11$ & $000$ & $001$ & $\cdots$
		\\ %\hline
		$I_{\{0,1\}^{\ast}\to\mathbb{Z}}(x)$ 
		& $0$ &$1$ & $-1$ &$2$ & $-2$
		& $3$ &$-3$& $4$ & $-4$ &  $\cdots$\\
		\hline
	\end{tabular}
%	
%	\vspace{5mm}
%	
%	\begin{tabular}{c|cccccccccc}
%		\hline
%		$x$ & $\lambda$ & $0$ & $1$ & $00$ & $01$ 
%		& $10$ & $11$ & $000$ & $001$ & $\cdots$
%		\\ \hline
%		$I_{\{0,1\}^{\ast}\to\mathbb{Z}}(x)$ 
%		& $0$ &$1$ & $-1$ &$2$ & $-2$
%		& $3$ &$-3$& $4$ & $-4$ &  $\cdots$\\
%		\hline
%	\end{tabular}
	\end{center}
For convenience, we define 
$I_{\sharp \to \mathcal{A}^{\ast}}:=I_{\mathcal{A}^{\ast}\to\sharp}^{-1}$.

The map $\mathcal{A}^{\ast}\times \mathcal{A}^{\ast}\ni (x,y)
\mapsto xy \in \mathcal{A}^{\ast}$
is called the \textit{concatenation}.
The set $\mathcal{A}^{\ast}$ with the concatenation is a monoid with 
identity element $\lambda$, i.e., 
$(xy)z=x(yz)$ for all $x,y,z\in \mathcal{A}^{\ast}$ and 
$\lambda x=x\lambda =x$ for all $x\in \mathcal{A}^{\ast}$.

The \textit{length} of $x\in \mathcal{A}^{\ast}$ is denoted by
$l(x)$ which is defined by 
$l(x)=n \stackrel{\mathrm{def}}{\Leftrightarrow} x\in \mathcal{A}^n$.
Obviously, we have for all $x, y \in \mathcal{A}^{\ast}, 
l(xy)=l(x)+l(y)$. 

For all $x,y \in \mathcal{A}^{\ast}$, we call $x$ a \textit{prefix} of $y$ 
if there exists  $z\in \mathcal{A}^{\ast}$ such that 
$y=xz$. 
A set $A \subset \mathcal{A}^{\ast}$ is said to be \textit{prefix-free} if, 
for all $x\in A$, the elements of $A\setminus \{x\}$ is not a prefix of $x$.
We set for all $x\in \mathcal{A}^{\ast}$
$$
	\Bar{x}:=\underbrace{1\cdots 1}_{l(x)}0x
$$
then we have $l(\bar{x})=2l(x)+1$.

Let $\mathcal{A}_1, \mathcal{A}_2$ be a nonempty finite set.
Let $\mathcal{D}$ be a subset of $\mathcal{A}_1^{\ast}$ 
and let $f$ be a function from $\mathcal{D}$ to $\mathcal{A}_2^{\ast}$.
If $\mathcal{D}\subsetneq \mathcal{A}_1^{\ast}$, 
we call such a function $f$ a 
\textit{partial function} and write 
$f:\mathcal{A}_1^{\ast}\rightsquigarrow\mathcal{A}_2^{\ast}$,
and if $\mathcal{D}= \mathcal{A}_1^{\ast}$ then
we call $f$ a \textit{total function}. 

A partial function 
$\phi :\mathcal{A}^{\ast}\rightsquigarrow\mathcal{A}^{\ast}$ is said to be 
\textit{partial recursive} if and only if there exists a Turing machine 
$M$ such that $\phi$ is computed by $M$, i.e., 
for all $x\in \mathcal{A}^{\ast}$, $M$ halts if and only if 
$x\in \mathrm{dom}(\phi)$, in that case, $M$ outputs $\phi(x)$. 

A partial function 
$\phi :\mathcal{A}_1^{\ast}\rightsquigarrow\mathcal{A}_2^{\ast}$ is partial 
recursive if there exists a partial recursive function 
$\psi :\mathcal{A}_1^{\ast}\rightsquigarrow\mathcal{A}_1^{\ast}$ such that 
$\phi = I_{\mathbb{Z}_+\to \mathcal{A}_2^{\ast}}\circ I_{\mathcal{A}_1^{\ast}\to\mathbb{Z}_+}\circ\psi$.

%A partial recursive function 
%$\phi :\mathcal{A}^{\ast}\rightsquigarrow\mathcal{A}^{\ast}$ 
%is called \textit{universal} if 
%there exists a universal Turing machine $U$ such that 
%$U$ computes $\phi$.

A \textit{partial recursive prefix function} 
$\phi :\mathcal{A}_1^{\ast}\rightsquigarrow\mathcal{A}_2^{\ast}$ is a 
partial recursive function such that $\mathrm{dom}(\phi)$ is prefix-free.

Let $\phi :\{0,1\}^{\ast}\rightsquigarrow\mathcal{A}^{\ast}$ be a 
partial recursive prefix function. For all $x\in \mathcal{A}^{\ast}$, 
the \textit{complexity} of $x$ with respect to $\phi$ is defined by 
$$
	K_{\phi}(x):=\begin{cases}
					\min\{l(p) : p\in \phi^{-1}(x)\}, & 
					(\phi^{-1}(x)\neq \emptyset),\\
					\infty & (\phi^{-1}(x)= \emptyset).
				\end{cases}
$$

A partial recursive prefix function 
$\phi :\{0,1\}^{\ast}\rightsquigarrow\mathcal{A}^{\ast}$ is said to be 
\textit{additively optimal} if for all partial recursive prefix function 
$\psi:\{0,1\}^{\ast}\rightsquigarrow\mathcal{A}^{\ast}$, 
there exists a constant $c_{\phi, \psi}\in \mathbb{R}$ such that 
$$
	\forall x\in \mathcal{A}^{\ast}, \quad
	K_{\phi}(x)\leq K_{\psi}(x)+c_{\phi,\psi}.
$$

\begin{theorem}
	There exists an additively optimal %universal 
	partial recursive prefix function.
\end{theorem}
{\it Proof}.
See \cite{LiVi}.
%Li Vitanyi Theorem 3.1.1.
\hfill $\square$

For each pair $(\phi, \psi)$ of additively optimal partial recursive prefix functions from 
$\{0,1\}^{\ast}$ to $\mathcal{A}^{\ast}$, there exists a constant 
$c_{\phi,\psi}>0$ such that for all $x\in \mathcal{A}^{\ast}$, 
$|K_{\phi}(x)-K_{\psi}(x)|\leq c_{\phi,\psi}$.

It is easily seen that any additively optimal partial recursive prefix function is surjective.

\begin{definition}
	We fix one additively optimal partial recursive prefix function 
	$\phi :\{0,1\}^{\ast}\rightsquigarrow\mathcal{A}^{\ast}$. 
	We define the prefix Kolmogorov complexity of $x\in \mathcal{A}^{\ast}$ 
	by 
	$$
		K(x):= K_{\phi}(x).
	$$
\end{definition}

%%%%%%%%%%%%%%%%%

\subsection{Shift dynamical system}

Let $\Sigma$ be a nonempty finite set, and we set $\Omega :=\Sigma^G$. 
By Tychonoff's theorem, $\Omega$ endowed with the product topology of 
the discrete topology on $\Sigma$ is a compact topological space.
It is well-known that this topology is also generated by the metric 
$d(\omega, \omega')=2^{-n(\omega, \omega')}$ where 
$n(\omega, \omega')=\sup\{n\in \mathbb{N} : \forall g\in \Lambda_n, 
\omega_g =\omega'_g\}$ for all $\omega=(\omega_g)_{g\in G}, 
\omega'=(\omega'_g)_{g\in G} \in \Omega$. 
For all $n\in \mathbb{N}$ and for all $s\in \Sigma^{\Lambda_n}$, we 
define the cylinder set of $s$ by 
$[\![s]\!]:=\{\omega\in \Omega : \omega \upharpoonright \Lambda_n=s\}$ 
where $\omega \upharpoonright \Lambda_n$ denotes the restriction of 
$\omega$ to $\Lambda_n$. 
Note that $[\![s]\!]$ is a clopen set. 
For all $n\in \mathbb{N}$, let $\mathcal{C}_n$ be the family of cylinder sets 
on $\Sigma^{\Lambda_n}$, i.e., 
$$
	\mathcal{C}_n:=\{[\![s]\!] : s \in \Sigma^{\Lambda_n}\},
$$
and set $\mathcal{C}:= \bigcup_n\mathcal{C}_n$.
The set 
$\mathcal{C}$ generates the Borel $\sigma$-algebra $\mathfrak{B}(\Omega)$.

We set 
$$
	\Sigma^{\Lambda_{\ast}}:=\bigcup_{n=0}^{\infty}\Sigma^{\Lambda_n}
$$
where $\Sigma^{\Lambda_0}:=\{\lambda\}$ and for all $n\in \mathbb{N}, 
\Sigma^{\Lambda_n}:=\{(\omega_g)_{g\in \Lambda_n} : \forall g\in \Lambda_n, 
\omega_g\in \Sigma\}$, and write 
$[\![V]\!]:=\bigcup_{s\in V}[\![s]\!]$ for all 
$V\subset \Sigma^{\Lambda_{\ast}}$.

Let $\sigma^g:\Omega \to \Omega$ denote the shift by $g\in G$, i.e., 
$(\sigma^g\omega)_i:=\omega_{i+g}$ for all $\omega=(\omega_g)_{g\in G}$, 
and we write $\sigma:=(\sigma^g)_{g\in G}$. 
Then $\sigma$ is a continuous action of $G$ on $\Omega$. Hence 
$(\Omega, \sigma)$ is a t.d.s.. 
% and which is called \textit{full shift}.
Note that $\sigma$ is a map from $G\times \Omega$ to $\Omega$, i.e.,
$\sigma:G\times \Omega \ni (g,\omega)\mapsto \sigma^g(\omega)\in \Omega$. 

A nonempty subset $S\subset \Omega$ is called a \textit{subshift} 
if and only if $S$ is shift-invariant (i.e. 
$\forall g \in G, \sigma^g(S)=S$) and $S$ is closed. 
If $S\subset \Omega$ is a subshift, then 
$(S, \sigma\upharpoonright (G\times S))$ is a t.d.s.. 
There exists a measure of maximal entropy measure for $\sigma\upharpoonright (G\times S)$
(see \cite{Ke}).

We call $f:G\to \mathbb{Z}_+$ a \textit{computable function} if 
there exists a partial recursive prefix function $\phi : \{0,1\}^{\ast}
\to \{0,1\}^{\ast}$ such that for all $(x_1,\cdots , x_d)\in G$, 
$$
	f(x_1,\cdots,x_d)=
	(I_{\{0,1\}^{\ast}\to\mathbb{Z}_+} \circ \phi)
	\left(\overline{I_{\sharp\to \{0,1\}^{\ast}}(x_1)}\cdots 
	\overline{I_{\sharp\to \{0,1\}^{\ast}}(x_{d-1})}
	I_{\sharp\to \{0,1\}^{\ast}}(x_d)
	\right)
$$
where 
$$
	\sharp =\begin{cases}
				\mathbb{Z}, & G=\mathbb{Z}^d, \\
				\mathbb{Z}_+, & G=\mathbb{Z}_+^d.
			\end{cases}
$$

We fix an arbitrary bijective computable function
$f:G\to \mathbb{Z}_+$ such that for all $n\in \mathbb{N}$, 
$$
	f(\Lambda_n)=\{0,1,\cdots, |\Lambda_n|-1\}
$$
and define $\mathcal{G}:\Sigma^{\Lambda_{\ast}}\to \Sigma^{\ast}$ as follows.
$$
	\mathcal{G}(s)
	:=\begin{cases}
		s_{f^{-1}(0)}\cdots s_{f^{-1}(|\Lambda_n|-1)}, & 
		s=(s_g)_{g\in\Lambda_n}\in \Sigma^{\Lambda_n}\ (n\in \mathbb{N}), \\
		\lambda, & s=\lambda.
	\end{cases}
$$

We define the prefix Kolmogorov 
complexity of $s\in \Sigma^{\Lambda_{\ast}}$ by 
$$
	\mathsf{K}(s):=K(\mathcal{G}(s)).
$$

\begin{lemma}\label{D_n}
	For all $n\in \mathbb{N}$ and $k\in \mathbb{R}_{\geq 0}$, we define 
	$$
		D_{n,k}:=\{s\in \Sigma^{\Lambda_n} : 
		\mathsf{K}(s)<k
		\}.
	$$
	Then %the following holds.
	$$
		|D_{n,k}|\leq 2^{k+1}.
	$$
\end{lemma}
{\it Proof.}
	Note that for all $s\in \Sigma^{\Lambda_n} \ (n\in \mathbb{N})$, 
	$\mathsf{K}(s)=K_{\phi}(\mathcal{G}(s))\neq \infty \ $. We define
	$\psi :\Sigma^{\Lambda_n}\to \{0,1\}^{\ast}$ such that 
	the following condition 
	holds.
	$$
		\psi(s)=p_s \Longleftrightarrow \phi(p_s)=\mathcal{G}(s) \ \text{and} 
		\ \mathsf{K}(s)
		=l(p_s).
	$$
	For all $s,t\in \Sigma^{\Lambda_n}$, we have
	$$
		s\neq t \Longrightarrow \phi(p_s)\neq\phi(p_t)
		\Longrightarrow p_s\neq p_t
	$$
	then $\psi$ is injective. Therefore 
	\begin{eqnarray*}
		|D_{n,k}|&=&|\{s\in \Sigma^{\Lambda_n} : \exists p_s \in \{0,1\}^{\ast}, 
		\psi(s)=p_s, l(p_s)<k\}|\\
		&\leq&|\{p\in\{0,1\}^{\ast} : l(p)<k\}| \\
		&\leq& 1+2+\cdots +2^{\lfloor k\rfloor}
		=2^{\lfloor k \rfloor +1}-1 \leq 2^{k+1}.
	\end{eqnarray*}
	\hfill $\square$
	
\begin{definition}[Kolmogorov complexity density]
	The upper and lower Kolmogorov complexity density 
	of $\omega \in \Omega$ are defined by
	$$
		\overline{\mathcal{K}}(\omega):=\limsup_{n\to \infty}
		\frac{\mathsf{K}(\omega \upharpoonright \Lambda_n)}{|\Lambda_n|},
		\quad
		\underline{\mathcal{K}}(\omega):=\liminf_{n\to \infty}
		\frac{\mathsf{K}(\omega \upharpoonright \Lambda_n)}{|\Lambda_n|}.
	$$
	If $\overline{\mathcal{K}}(\omega)=\underline{\mathcal{K}}(\omega)$, 
	we simply denote them 
	by $\mathcal{K}(\omega)$, i.e.,
	$$
		\mathcal{K}(\omega):=\lim_{n\to \infty}
		\frac{\mathsf{K}(\omega \upharpoonright \Lambda_n)}{|\Lambda_n|}.
	$$
	\end{definition}
	
\begin{remark}
	The quantities $\overline{\mathcal{K}}(\omega)$ and 
	$\underline{\mathcal{K}}(\omega)$ are independent of the choice of 
	additively optimal partial recursive prefix function $\phi$ and 
	$\mathcal{G}$ and uniquely defined.
\end{remark}

%%%%%%%%%%%%%%%%%%%

\section{Relation between KS entropy and Kolmogorov complexity}

%In this section, we prove the main theorem. 
Let $d\in \mathbb{N}$, $G=\mathbb{Z}^d$ or $G=\mathbb{Z}^d_+$, 
$\Sigma$ be a nonempty finite set, and $S\subset \Omega \ (:=\Sigma^G)$ be a subshift.
Other notations are the same as before. 
We set $\varsigma:= \sigma\upharpoonright (G\times S)$. 
%i.e., 
%
%and let
%$S\subset \Omega$ be a subshift. 
Note that $(S, \varsigma)$ is a t.d.s. 
We now state the main result.

%Before prove the main theorem, we introduce Brudno's theorem for shift space.
%
%\begin{theorem}[Brudno's theorem]
%	Let $G=$
%	If $\mu\in EM(\Omega, \sigma)$, then $\overline{\mathcal{K}}(\omega)=h_{\sigma}(\mu)$
%	for $\mu$-a.e. $\omega\in\Omega$.
%\end{theorem}
%{\it Proof}. 
%See \cite{Br}.

\begin{theorem}[Brudno's theorem for $\mathbb{Z}^d$ (or 
	$\mathbb{Z}^d_+$) subshifts]\label{main}
	If $\mu \in EM(S, \varsigma)$, then 
	\begin{equation}
		%\lim_{n\to \infty}
		%\frac{\mathsf{K}(\omega \upharpoonright \Lambda_n)}{|\Lambda_n|}
		\mathcal{K}(\omega)=h_{\varsigma}(\mu), \quad
		\mu\text{-a.e.} \omega \in S.
	\end{equation}\label{Brudno}
\end{theorem}

\begin{remark}
	Brudno's original result is on the case $G=\mathbb{Z}_+$ only 
	\cite{Br}. 
	In the case $G=\mathbb{Z}^d$ or $G=\mathbb{Z}^d_+$, Simpson showed that if 
	%The special case that 
	$\mu$ is a measure of maximal entropy, 
	then (\ref{Brudno}) holds \cite{Sim}.
	%is shown by Simpson \cite{Sim}.
	Our theorem is a generalization of them.
\end{remark}

We prove Theorem \ref{main} by giving two lemmas.

\begin{lemma}\label{infmain}
	If $\mu \in EM(S, \varsigma)$, then 
	\begin{equation}
		\underline{\mathcal{K}}(\omega)\geq h_{\varsigma}(\mu), \quad
		\mu\text{-a.e.} \omega \in S.
		\label{infBrudno}
	\end{equation}
\end{lemma}
{\it Proof.}
If $h_{\varsigma}(\mu)=0$, then (\ref{infBrudno}) is obvious.

Let $h_{\varsigma}(\mu)>0$ and fix an arbitrary $k\in \mathbb{N}$ 
such that $\frac{1}{k}<h_{\varsigma}(\mu)$.
For all $n \in \mathbb{N}$, we set
$$
	\widetilde{{D}_{n,k}}:=\left\{
	s\in\Sigma^{\Lambda_n} : \frac{\mathsf{K}(s)}{|\Lambda_n|}\leq
	h_{\varsigma}(\mu)-\frac{1}{k}
	\right\}.
$$
By Lemma \ref{D_n}, we have 
\begin{equation}
|\widetilde{D_{n,k}}|\leq 2^{|\Lambda_n|(h_{\varsigma}(\mu)-\frac{1}{k})+1}.
\label{|D_n|}
\end{equation}
We fix an arbitrary $\epsilon \in \left(0,\frac{1}{k}\right)$ and set
$$
	T_{n,k,\epsilon}:=\left\{
					s\in \Sigma^{\Lambda_n} : \mu([\![s]\!] \cap S)
					< 2^{-|\Lambda_n|(h_{\varsigma}(\mu)-\frac{1}{k}+
					\epsilon)}
					\right\}.
$$
By Shannon-McMillan-Breiman theorem (Theorem \ref{SMB}), 
the following holds for $\mu$-a.e. $\omega \in S$
$$
	\exists N_{\omega}, \forall n\geq N_{\omega}, 
	\left|
		h_{\varsigma}(\mu)-\frac{-\log_2\mu([\![\omega\upharpoonright 
		\Lambda_n]\!] \cap S)}{|\Lambda_n|}
	\right|
	<\frac{1}{k}-\epsilon.
$$
Note that
\begin{eqnarray*}
	\left|
		h_{\varsigma}(\mu)-\frac{-\log_2\mu([\![\omega\upharpoonright 
		\Lambda_n]\!] \cap S)}{|\Lambda_n|}
	\right|
	<\frac{1}{k}-\epsilon
	&\Longrightarrow &
	\frac{-\log_2\mu([\![\omega\upharpoonright 
		\Lambda_n]\!] \cap S)}{|\Lambda_n|}
		>h_{\varsigma}(\mu)-\frac{1}{k}+\epsilon \\
	&\Longleftrightarrow&
	2^{-|\Lambda_n |(h_{\varsigma}(\mu)-\frac{1}{k}+\epsilon)}
	>\mu([\![\omega\upharpoonright 
		\Lambda_n]\!] \cap S) \\
	&\Longleftrightarrow&
	\omega \upharpoonright \Lambda_n \in T_{n,k,\epsilon} \\
	&\Longleftrightarrow&
	\omega \in [\![T_{n,k,\epsilon}]\!]\cap S.
\end{eqnarray*}
Hence we have for $\mu$-a.e. $\omega \in S$
\begin{equation}
	\exists N_{\omega}, \forall n\geq N_{\omega}, 
	\omega \in [\![T_{n,k,\epsilon}]\!]\cap S.
	\label{Tnke}
\end{equation}
On the other hand, by (\ref{|D_n|}) and the definition of $T_{n,k,\epsilon}$, 
we have
\begin{eqnarray*}
	\mu([\![\widetilde{D_{n,k}}]\!]\cap [\![T_{n,k,\epsilon}]\!]\cap S)
	%=\mu([\![D_{n,k}\cap T_{n,k,\epsilon}]\!]\cap S)
	&=&\mu\left(
	\bigcup_{s\in \widetilde{D_{n,k}}\cap T_{n,k,\epsilon}}[\![s]\!]\cap S
	\right)\\
	&\leq&
	\sum_{s\in \widetilde{D_{n,k}}\cap T_{n,k,\epsilon}}\mu([\![s]\!]\cap S)\\
	&\leq&
	2^{|\Lambda_n|(h_{\varsigma}(\mu)-\frac{1}{k})+1}\cdot
	2^{-|\Lambda_n|(h_{\varsigma}(\mu)-\frac{1}{k}+
					\epsilon)}
					=2^{-|\Lambda_n|\epsilon +1}.
\end{eqnarray*}
Hence
\begin{eqnarray*}
	\sum_{n=1}^{\infty}\mu([\![\widetilde{D_{n,k}}]\!]\cap [\![T_{n,k,\epsilon}]\!]\cap S)
	<\infty.
\end{eqnarray*}
Therefore, by the Borel-Cantelli lemma,
for $\mu$-a.e. $\omega \in S$,
\begin{equation}
	\exists N'_{\omega}\in\mathbb{N}, \forall n
	\geq N'_{\omega}, \omega \notin [\![\widetilde{D_{n,k}}]\!]\cap [\![T_{n,k,\epsilon}]\!]\cap S.
	\label{BC}
\end{equation}
By (\ref{Tnke}) and (\ref{BC}), for $\mu$-a.e. $\omega \in S$, we have
\begin{equation}
	\exists N''_{\omega}\in\mathbb{N}, \forall n
	\geq N''_{\omega}, \omega \notin [\![\widetilde{D_{n,k}}]\!].
\end{equation}
Since $\omega \notin [\![\widetilde{D_{n,k}}]\!]$ means 
$\frac{\mathsf{K}(\omega\upharpoonright \Lambda_n)}{|\Lambda_n|}>
	h_{\varsigma}(\mu)-\frac{1}{k}$, for all $k\in \mathbb{N}$, we 
	have 
		\begin{equation}
		\underline{\mathcal{K}}(\omega)\geq h_{\varsigma}(\mu)
		-\frac{1}{k}, \quad
		\mu\text{-a.e.} \omega \in S.
	\end{equation}
Thus (\ref{infBrudno}) holds.
\hfill $\square$

\begin{lemma}\label{supmain}
	If $\mu \in EM(S, \varsigma)$, then 
	\begin{equation}
		\overline{\mathcal{K}}(\omega)\leq h_{\varsigma}(\mu), \quad
		\mu\text{-a.e.} \omega \in S.
		\label{supBrudno}
	\end{equation}
\end{lemma}
{\it Proof.}
Fix an arbitrary $m\in \mathbb{N}$ and let $L_m$ be the side length of the 
hypercube $\Lambda_m$, i.e.,
$$
	L_m:=\begin{cases}
		m & \text{if}\ G=\mathbb{Z}_+^d, \\
		2m-1 & \text{if}\ G=\mathbb{Z}^d.
		\end{cases}
$$
For each $n\in \mathbb{N}_{> m}$, let us consider a covering of $\Lambda_n$ by 
shifted $\Lambda_m$. In particular, there uniquely exists $k\in \mathbb{N}$ 
such that 
$$
	\bigsqcup_{g\in \Lambda_k}(L_m g+\Lambda_m) \subsetneq
	\Lambda_n \subset 
	\bigsqcup_{g\in \Lambda_{k+1}}(L_m g+\Lambda_m),
$$
where $L_mg+\Lambda_m:=\{L_mg+h : h\in \Lambda_m\}$ and $\bigsqcup$ denotes 
disjoint union.
We set $\check{\Lambda}_n:=\bigsqcup_{g\in \Lambda_k}(L_m g+\Lambda_m)$ and  
$M:=|\Sigma^{\Lambda_m}|$.
Let us consider bijective map $r:\{1,\cdots, |\Sigma^{\Lambda_m}|\}\to 
\Sigma^{\Lambda_m}$, and let $r_j:=r(j)\ (1\leq j\leq M)$. 
For an arbitrary $\omega \in S$, we define 
\begin{eqnarray*}
	\mathsf{f}_{r_j}(\omega)&:=&|\{g\in \Lambda_k : \varsigma^{-L_mg}\omega
	\in [\![r_j]\!]\cap S\}|,\\
	\mathsf{f}_r(\omega)&:=&(\mathsf{f}_{r_1}(\omega),\cdots, 
	\mathsf{f}_{r_M}(\omega))
	\in \mathbb{Z}_+^M.
\end{eqnarray*}
By definition, $\mathsf{f}_{r_1}(\omega)+\cdots +
	\mathsf{f}_{r_M}(\omega)=|\Lambda_k|$
	holds for all $\omega \in S$. 
Since for all $s\in \Sigma^{\check{\Lambda}_n}$ and for all 
$\omega_1,\omega_2 \in [\![s]\!]\cap S$, we have 
$\mathsf{f}_{r_j}(\omega_1)=\mathsf{f}_{r_j}(\omega_2)$. 
We set $\mathsf{f}_{r_j}(s):=\mathsf{f}_{r_j}(\omega)\ (\omega \in [\![s]\!]\cap S)$.
We endow $\Sigma^{\check{\Lambda}_n}$ with an equivalence relation as follows:
$$
	\forall s_1,s_2\in \Sigma^{\check{\Lambda}_n}, \quad 
	s_1 
	%\stackrel{\mathsf{f}_r}{\sim} 
	\sim_{\mathsf{f}_r} 
	s_2 
	\stackrel{\mathrm{def}}{\Longleftrightarrow} 
	\mathsf{f}_r(s_1)=\mathsf{f}_r(s_2).
$$
For $s\in \Sigma^{\check{\Lambda}_n}$, let 
$[s]_{\mathsf{f}_r}:=\{t\in 
\Sigma^{\check{\Lambda}_n} : s \sim_{\mathsf{f}_r} t
\}$ be an equivalence class 
of $s$ by $\sim_{\mathsf{f}_r} $.
Then we have
$$
	|[s]_{\mathsf{f}_r}|=\frac{|\Lambda_k|!}
	{\mathsf{f}_{r_1}(s)!\mathsf{f}_{r_2}(s)!\cdots \mathsf{f}_{r_M}(s)!}.
$$
By the above mentioned preparations, 
we take the following procedures:
%let us consider the followings:
% three bijective maps.
\begin{enumerate}
	\item We fix a bijective map $\mathcal{F}$ from 
	$\Sigma^{\check{\Lambda}_n}/ \sim_{\mathsf{f}_r}$ to 
	$V_k^M:=
	\{(x_1,\cdots,x_M)\in \mathbb{Z}_+^M: x_1+\cdots +x_M=|\Lambda_k|\}$
	such that
	$$
		\mathcal{F}:[s]_{\mathsf{f}_r}\mapsto 
(\mathsf{f}_{r_1}(s),\cdots, \mathsf{f}_{r_M}(s)),
	$$
	and identify $\Sigma^{\check{\Lambda}_n}/ \sim_{\mathsf{f}_r}$
	with $V_k^M$.
	\item We fix the arbitrary bijective maps 
	$\mathcal{N}_{[s]_{\mathsf{f}_r}}$ and 
	$\mathcal{R}_{\Lambda_n\setminus \check{\Lambda}_n}$
	such that 
	for each $[s]_{\mathsf{f}_r}\in \Sigma^{\check{\Lambda}_n}/
	\sim_{\mathsf{f}_r}$, 
	$$
		\mathcal{N}_{[s]_{\mathsf{f}_r}}:[s]_{\mathsf{f}_r}\to \{1,\cdots, |[s]_{\mathsf{f}_r}|\}
	$$ 
	and
	$$
		\mathcal{R}_{\Lambda_n\setminus \check{\Lambda}_n}:
		\Sigma^{\Lambda_n\setminus \check{\Lambda}_n}\to
		\{1,2,\cdots, |\Sigma^{\Lambda_n\setminus \check{\Lambda}_n}|\}.
	$$
%	\item fix the arbitrary bijective map
\end{enumerate}
Then %for each
we can uniquely identify each 
$t \in \Sigma^{\Lambda_n}$ 
%there uniquely exists
with
$(\mathbf{x}, y, z)$ 
%such that 
where $\mathbf{x}\in V_k^M$, 
$y\in \mathcal{N}_{\mathcal{F}^{-1}(\mathbf{x})}(\mathcal{F}^{-1}(\mathbf{x}))
$ and $z\in \{1,2,\cdots, |\Sigma^{\Lambda_n\setminus \check{\Lambda}_n}|\}$. 
Hence there exists a partial recursive prefix function 
$\phi:\{0,1\}^{\ast}\to \Sigma^{\ast}$ such that 
\begin{eqnarray*}
	\forall \omega \in S, \exists \mathbf{x}=(x_1,\cdots, x_M)\in V_k^M, 
	\exists y\in \mathcal{N}_{\mathcal{F}^{-1}(\mathbf{x})}(\mathcal{F}^{-1}(\mathbf{x})), \exists z\in \{1,2,\cdots, |\Sigma^{\Lambda_n\setminus \check{\Lambda}_n}|\},\\
	\phi\left(
	\overline{I_{\mathbb{Z}_+\to \{0,1\}^{\ast}}(x_1)}\cdots
	\overline{I_{\mathbb{Z}_+\to \{0,1\}^{\ast}}(x_M)}\ 
	\overline{I_{\mathbb{Z}_+\to \{0,1\}^{\ast}}(y)}\ 
	I_{\mathbb{Z}_+\to \{0,1\}^{\ast}}(z)
	\right)
	=\mathcal{G}(\omega \upharpoonright \Lambda_n). 
\end{eqnarray*}
Obviously, $x_j=\mathsf{f}_{r_j}(\omega)\ (j\in \{1,\cdots, M\})$.
By the definition of $K_{\phi}$, we have
\begin{eqnarray}
	K_{\phi}(\mathcal{G}(\omega\upharpoonright \Lambda_n))&\leq& 
	l\left(
	\overline{I_{\mathbb{Z}_+\to \{0,1\}^{\ast}}(x_1)}\cdots
	\overline{I_{\mathbb{Z}_+\to \{0,1\}^{\ast}}(x_M)}\ 
	\overline{I_{\mathbb{Z}_+\to \{0,1\}^{\ast}}(z)}\ 
	I_{\mathbb{Z}_+\to \{0,1\}^{\ast}}(y)
	\right)\notag \\
	&\leq&
	2\left(
	\sum_{j=1}^Ml(I_{\mathbb{Z}_+\to \{0,1\}^{\ast}}(x_j))
	+l(I_{\mathbb{Z}_+\to \{0,1\}^{\ast}}(z))
	\right)\notag\\
	&&+M+1+
	l(I_{\mathbb{Z}_+\to \{0,1\}^{\ast}}(y)).\label{Kp}
\end{eqnarray}
The following inequalities can be easily seen:
\begin{eqnarray*}
	&&\sum_{j=1}^Ml(I_{\mathbb{Z}_+\to \{0,1\}^{\ast}}(x_j))\leq 
	\sum_{j=1}^M\log_2(x_j+1)\leq
	M\log_2(|\Lambda_k|+1),\\
	&&l(I_{\mathbb{Z}_+\to \{0,1\}^{\ast}}(z))\leq
	\log_2(|\Sigma^{\Lambda_n\setminus\check{\Lambda}_n}|+1)\leq
	(|\Lambda_n|-|\check{\Lambda}_n|)\log_2|\Sigma|+1,\\
	&&l(I_{\mathbb{Z}_+\to \{0,1\}^{\ast}}(y))\leq 
	\log_2 y+1\leq
	\log_2\frac{|\Lambda_k|!}{x_1!\cdots x_M!}+1.
\end{eqnarray*}
By these inequalities and (\ref{Kp}), we have 
\begin{equation}
	\frac{K_{\phi}(\mathcal{G}(\omega\upharpoonright \Lambda_n))}
	{|\Lambda_n|}\leq 
	\frac{1}{|\Lambda_n|}\log_2\frac{|\Lambda_k|!}{x_1!\cdots x_M!}+
	\frac{2(|\Lambda_n|-|\check{\Lambda}_n|)}{|\Lambda_n|}\log_2|\Sigma|
	+o(1), \quad (n\to \infty).\label{Kp2}
\end{equation}
Let us estimate the right hand side of (\ref{Kp2}).
By direct computations using Stirling's formula, we can see that
\begin{equation}
	\frac{1}{|\Lambda_n|}\log_2\frac{|\Lambda_k|!}{x_1!\cdots x_M!}
	=\frac{|\Lambda_k|}{|\Lambda_n|}\sum_{j=1}^M\varphi\left(
	\frac{x_j}{|\Lambda_k|}\right)+o(1), \quad (n\to \infty).\label{St}
\end{equation}
Since $x_j=\mathsf{f}_{r_j}(\omega)$ and 
$\mu \in EM(S,\varsigma)$,
%the definition of identification $\mathcal{F}$ and 
%$\mathsf{f}_{r_j}$, 
for $\mu$-a.e.$\omega\in S$,
we have the following.
\begin{eqnarray}
	\lim_{k\to\infty}\frac{x_j}{|\Lambda_k|}&=&
	\lim_{k\to\infty}\frac{1}{|\Lambda_k|}
	|\{g\in \Lambda_k : \varsigma^{-L_mg}\omega \in [\![r_j]\!]\cap S\}|\notag
	\\
	&=&
	\lim_{k\to\infty}\frac{1}{|\Lambda_k|}\sum_{g\in\Lambda_k}1_{[\![r_j]\!]\cap S}
	(\varsigma^{-L_mg}\omega)=
	\mu([\![r_j]\!]\cap S).\label{001}
\end{eqnarray}
The last equality is derived from the Birkhoff's ergodic theorem.

By $\Lambda_{n-L_m}\subset \check{\Lambda}_n\subsetneq \Lambda_n$ and 
$|\check{\Lambda}_n|=|\Lambda_m|\cdot |\Lambda_k|$, we have 
$|\Lambda_{n-L_m}|\leq |\Lambda_m|\cdot |\Lambda_k|\leq |\Lambda_n|$.
Then 
$$
	\frac{1}{|\Lambda_m|}\cdot \frac{|\Lambda_{n-L_m}|}{|\Lambda_n|}
	\leq \frac{|\Lambda_k|}{|\Lambda_n|}
	\leq \frac{1}{|\Lambda_m|}
$$
and we have
$$
	\frac{|\Lambda_{n-L_m}|}{|\Lambda_n|}
	=\begin{cases}
		\frac{(n-L_m)^d}{n^d}\to 1 & \text{if} \ G=\mathbb{Z}_+^d,\\
		\frac{\{2(n-L_m)-1\}^d}{(2n-1)^d} \to 1 & \text{if} \ G=\mathbb{Z}^d,
	\end{cases}
	\quad (n\to \infty).
$$
Hence the following hold:
\begin{equation}
	\lim_{n\to\infty}\frac{|\Lambda_k|}{|\Lambda_n|}=\frac{1}{|\Lambda_m|}, \
	\lim_{n\to\infty}\frac{|\check{\Lambda}_n|}{|\Lambda_n|}=1.\label{002}
	\end{equation}
Obviously, if $n\to \infty$, then $k\to\infty$.
By (\ref{Kp2}), (\ref{St}), (\ref{001}), (\ref{002}) and 
Kolmogorov complexity's definition, for $\mu$-a.e.$\omega \in S$ and for 
all $m\in \mathbb{N}$, we have
\begin{eqnarray*}
	\overline{\mathcal{K}}(\omega)&\leq& 
	\limsup_{n\to\infty}\frac{K_{\phi}
	(\mathcal{G}(\omega\upharpoonright \Lambda_n))}{|\Lambda_n|}\\
	&\leq&
	\limsup_{n\to\infty}
	\left\{
	\frac{|\Lambda_k|}{|\Lambda_n|}\sum_{j=1}^M\varphi\left(
	\frac{x_j}{|\Lambda_k|}\right)+
	\frac{2(|\Lambda_n|-|\check{\Lambda}_n|)}{|\Lambda_n|}\log_2|\Sigma|
	\right\}\\
	&=&
	\frac{1}{|\Lambda_m|}\sum_{j=1}^M\varphi(\mu([\![r_j]\!]\cap S)).
\end{eqnarray*}
Hence we have the following inequality:
\begin{equation}
	\overline{\mathcal{K}}(\omega)\leq
	\liminf_{m\to\infty}
	\frac{1}{|\Lambda_m|}\sum_{j=1}^M\varphi(\mu([\![r_j]\!]\cap S)), \quad
	\text{for}\ \mu\text{-a.e.}\ \omega \in S.\label{l1}
\end{equation}
Note that $\alpha:=\{[\![\omega \upharpoonright \Lambda_1]\!]\cap S\}
_{\omega \in S}$ 
is a $\mu$-generator and 
$$
	\alpha^{\Lambda_m}=\bigvee_{g\in \Lambda_m}\varsigma^{-g}\alpha 
	=\{[\![r_j]\!]\cap S\}_{j=1}^M.
$$
Therefore, by Kolmogorov-Sinai theorem, we have the following equation.
\begin{equation}
	\lim_{m\to\infty}
	\frac{1}{|\Lambda_m|}\sum_{j=1}^M\varphi(\mu([\![r_j]\!]\cap S))=
	\lim_{m\to\infty}\frac{1}{|\Lambda_m|}H_{\mu}(\alpha^{\Lambda_m})
	=h_{\varsigma}(\mu).\label{l2}
\end{equation}
(\ref{l1}) and (\ref{l2}) complete the proof.
\hfill $\square$

\medskip

Theorem \ref{main} follows from Lemma \ref{infmain} and Lemma \ref{supmain}.

\begin{example}[$d$-dimensional Bernoulli shifts]
Let $(\Omega, \sigma)$ be the $d$-dimensional shift space as before.
We fix a probability vector $q=(q_{i} : i \in \Sigma)$ on $\Sigma$ and 
denote the corresponding Bernoulli measure on $\mathfrak{B}(\Omega)$ by $\mu:=q^{\times G}$.
Then, by Kolmogorov-Sinai theorem, we can show that
$
	h_{\sigma}(\mu)=\sum_{i\in \Sigma}\varphi(q_i).
$
By Theorem \ref{main}, we have for $\mu$-a.e. $\omega\in \Omega$ 
$$
	\mathcal{K}(\omega)=\sum_{i\in \Sigma}\varphi(q_i).
$$
\end{example}

\section{Representation of the pressure by using Kolmogorov complexity density}

Let us consider some applications of Theorem \ref{main}. 
%Although this section's result are known in the case $d=1$, 
%but by 
%Several interesting fact is found by its generalized theorem.
Notations are the same as in Section 3.

\begin{theorem}\label{gBrudno}
	If $\mu \in M(S, \varsigma)$, then we have
	\begin{equation}
		h_{\varsigma}(\mu)=
		\mu(\mathcal{K})
		%\int_{S}\mathcal{K}(\omega)d\mu (\omega)
		=\lim_{n\to \infty}\frac{1}{|\Lambda_n|}\sum_{s\in \Lambda_n}
		\mathsf{K}(s)\mu ([\![s]\!]\cap S).
	\end{equation}
\end{theorem}
{\it Proof}.
Let $\mu=\int_{EM(S,\varsigma)}\nu d\rho(\nu)$ be the ergodic decomposition. 
By Theorem \ref{ergdecomp}, Jacobs's theorem (Theorem \ref{Jacob}) and 
the main theorem (Theorem \ref{main}), 
we have 
	$$
		\int_{S}\mathcal{K}(\omega)d\mu (\omega)
		=\int_{EM(S,\varsigma)}\left\{
		\int_S \mathcal{K}(\omega)d\nu(\omega)
		\right\}
		d\rho(\nu)
		=\int_{EM(S,\varsigma)}
		h_{\varsigma}(\nu)
		d\rho(\nu)
		=h_{\varsigma}(\mu).
	$$
On the other hand, by Lebesgue's convergence theorem, we have
	\begin{eqnarray*}
		\lim_{n\to \infty}\frac{1}{|\Lambda_n|}\sum_{s\in \Lambda_n}
		\mathsf{K}(s)\mu ([\![s]\!]\cap S)
		&=&\lim_{n\to \infty}\int_S\frac{\mathsf{K}(\omega\upharpoonright \Lambda_n)}{|\Lambda_n|}
		d\mu(\omega)\\
		&=&\int_S\lim_{n\to \infty}\frac{\mathsf{K}(\omega\upharpoonright \Lambda_n)}{|\Lambda_n|}
		d\mu(\omega)
		=\int_{S}\mathcal{K}(\omega)d\mu (\omega).
	\end{eqnarray*}
\hfill $\square$

\medskip

\begin{remark}
	In the case $G=\mathbb{Z}_+$, 
	Theorem \ref{gBrudno} can be found in \cite{BBG}.
\end{remark}

Theorem \ref{gBrudno} immediately leads one to the following theorem.

\begin{theorem}[Variational principle]\label{pressure}
		Let $\psi \in USC(S), \ \inf \psi >-\infty$.
	Then the pressure of $\psi$ is given by
		$$
			p(\psi)=\sup_{\mu\in M(S, \varsigma)}\mu(\mathcal{K}+\psi).
		$$
	In particular, the topological entropy is $\sup_{\mu\in M(S, \varsigma)}\mu(\mathcal{K})$. 
%	If $\mu\in M(S,\varsigma)$ is an equilibrium state for $\psi \in USC(S)$ then
%		$$
%			p(0)=h_{\varsigma}(\mu)=\mu(\mathcal{K}).
%		$$
	If $\mu\in M(S,\varsigma)$ is an equilibrium state for $\psi$, 
	then we have
	$$
		p(\psi)=\mu(\mathcal{K}+\psi).
	$$
\end{theorem}

Theorem \ref{pressure} shows that, in an equilibrium state, 
the pressure means the expectation value of the sum of 
Kolmogorov complexity density and local energy.
%Obviously, 
For example, this theorem is directly 
applicable to the $d$-dimensional Ising model.

%By using Theorem \ref{pressure}, 
%we can express the pressure of the $d$-dimensional Ising model by 
%using Kolmogorov complexity density.

%\begin{example}[Ising model]
%	Let $d$ be a positive integer and let
%	$\Sigma:=\{-1,+1\}$ represent ``spin up'' and ``spin down'' at each lattice point of 
%	$G:=\mathbb{Z}^d$. We set $\Omega:=\Sigma^{G}$ and let $\mathcal{T}$ be a shift action of $G$ 
%	on $\Omega$. We define local energy function $\psi_{\beta, B}:\Omega \to \mathbb{R}$ 
%	such that for $\omega\in \Omega$
%	$$
%		\psi_{\beta, B}(\omega):=-\beta(\psi(\omega)-B\omega_0)
%	$$
%	where $\psi(\omega):=-\sum_{j=1}^d(\omega_0\omega_{e_j}+\omega_0\omega_{-e_j})$ 
%	denote the interaction neighboring spins, 
%	$-B\omega_0$ denote the effect of a magnetic field $B\in \mathbb{R}$ on the spin at site $0$
%	and $\beta\geq 0$ denote the inverse temperature. 
%	Note that there exists a equilibrium state (see \cite{Ke}). 
%	Let $\mu$ be a equilibrium state, then the pressure of this model is given by
%	$$
%		p(\psi)=\mu(\mathcal{K}+\psi).
%	$$
%	
%\end{example}

%Then we have the injective map such that for all $s \in \Sigma^{\Lambda_n}$, 
%$$
%	\mathcal{G}(s):=(\mathcal{F}([s\upharpoonright 
%	\Sigma^{\check{\Lambda}_n}]_{\mathsf{f}_r}), 
%	N_{[s\upharpoonright \Sigma^{\check{\Lambda}_n}]_{\mathsf{f}_r}}
%	(s\upharpoonright \Sigma^{\check{\Lambda}_n}), 
%	\mathcal{R}_{\Lambda_n\setminus \check{\Lambda}_n}
%	(s\upharpoonright \Sigma^{\Lambda_n\setminus \check{\Lambda}_n}))
%$$

%from $\Sigma^{\Lambda_n}$ to 
%$V_k^M$
%$\{(x_1,\cdots,x_M)\in \mathbb{Z}_+^M: x_1+\cdots +x_M=|\Lambda_k|\}
%\times $

%%%%%%%%%%%%%

\section*{Acknowledgments}
We are grateful to A. Arai for reading our paper and for his valuable comments.
We also thank K. Tadaki for his important comments about how we should 
define the prefix Kolmogorov complexity of $s\in \Sigma^{\Lambda_{\ast}}$.
We would also like to express our sincere gratitude to M. Yuri for her 
valuable comments and longstanding encouragements. 
T. Fuda would like to thank S. Galatolo and A. Sakai for valuable discussions.
%We thank A. Arai, S. Galatolo, A. Sakai, K. Tadaki and M. Yuri for valuable comments.
%The author would like to thank 
%Professor Asao Arai for valuable comments.

\end{document}